\documentclass[12pt]{amsart}
\usepackage{amscd}
\def\NZQ{\Bbb}               
\def\NN{{\NZQ N}}

\def\frk{\frak}               

\def\mm{{\frk m}}

\def\Phi{{\frk n}}
\def\Phi{{\frk N}}
\def\opn#1#2{\def#1{\operatorname{#2}}} 
\opn\chara{char}
\opn\length{\ell}
\opn\pd{pd}
\opn\rk{rk}
\opn\projdim{proj\,dim}
\opn\injdim{inj\,dim}
\opn\rank{rank}
\opn\depth{depth}
\opn\grade{grade}
\opn\height{height}
\opn\embdim{emb\,dim}
\opn\codim{codim}

\opn\Tr{Tr}
\opn\bigrank{big\,rank}
\opn\superheight{superheight}\opn\lcm{lcm}
\opn\trdeg{tr\,deg}%
\opn\reg{reg}
\opn\lreg{lreg}
\opn\ini{in}
\opn\lpd{lpd}
\opn\div{div}
\opn\Div{Div}
\opn\cl{cl}
\opn\Cl{Cl}
\opn\Spec{Spec}
\opn\Supp{Supp}
\opn\supp{supp}
\opn\Sing{Sing}
\opn\Ass{Ass}
\opn\Ann{Ann}
\opn\Rad{Rad}
\opn\Soc{Soc}
\opn\Im{Im}
\opn\Ker{Ker}
\opn\Coker{Coker}
\opn\Am{Am}
\opn\Hom{Hom}
\opn\Tor{Tor}
\opn\Ext{Ext}
\opn\End{End}
\opn\Aut{Aut}
\opn\id{id}

\opn\nat{nat}
\opn\pff{pf}
\opn\Pf{Pf}
\opn\GL{GL}
\opn\SL{SL}
\opn\mod{mod}
\opn\ord{ord}
\opn\Gin{Gin}
\opn\aff{aff}
\opn\con{conv}
\opn\relint{relint}
\opn\st{st}
\opn\lk{lk}
\opn\cn{cn}
\opn\core{core}
\opn\vol{vol}
\opn\link{link}
\opn\star{star}
\opn\gr{gr}

\def\pot#1#2{#1[\kern-0.28ex[#2]\kern-0.28ex]}

\opn\dirlim{\underrightarrow{\lim}}
\opn\inivlim{\underleftarrow{\lim}}
\let\dirsum=\oplus

\let\Dirsum=\bigoplus

\let\to=\rightarrow
\let\To=\longrightarrow
\def\Implies{\ifmmode\Longrightarrow \else
       \unskip${}\Longrightarrow{}$\ignorespaces\fi}
\def\implies{\ifmmode\Rightarrow \else
       \unskip${}\Rightarrow{}$\ignorespaces\fi}
\def\iff{\ifmmode\Longleftrightarrow \else
       \unskip${}\Longleftrightarrow{}$\ignorespaces\fi}

\let\:=\colon
\newtheorem{Theorem}{Theorem}[section]
\newtheorem{Lemma}[Theorem]{Lemma}
\newtheorem{Corollary}[Theorem]{Corollary}
\newtheorem{Proposition}[Theorem]{Proposition}
\newtheorem{Remark}[Theorem]{Remark}

\newtheorem{Example}[Theorem]{Example}

\newtheorem{Definition}[Theorem]{Definition}

\let\epsilon\varepsilon
\let\phi=\varphi
\let\kappa=\varkappa
\def\qed{\ifhmode\textqed\fi
     \ifmmode\ifinner\quad\qedsymbol\else\dispqed\fi\fi}
\def\textqed{\unskip\nobreak\penalty50
      \hskip2em\hbox{}\nobreak\hfil\qedsymbol
      \parfillskip=0pt \finalhyphendemerits=0}
\def\dispqed{\rlap{\qquad\qedsymbol}}

\opn\dis{dis}
\def\pnt{{\raise0.5mm\hbox{\large\bf.}}}

\opn\Lex{Lex}

\begin{document}

\title{Rigid resolutions and big Betti numbers }

\author{Aldo Conca, J\"urgen Herzog and  Takayuki Hibi}
\address{Aldo Conca, DIMA - Dipartimento di Matematica,
Via Dodecaneso, 35, 16146 Genova, Italy}
\email{conca@dima.unige.it}
\address{J\"urgen Herzog, Fachbereich Mathematik und
Informatik, Universit\"at Duisburg-Essen, Campus Essen,
45117 Essen, Germany}
\email{juergen.herzog@uni-essen.de}
\address{Takayuki Hibi, Department of Mathematics, Graduate School of Science,
Osaka University, Toyonaka, Osaka 560-0043, Japan}
\email{hibi@math.sci.osaka-u.ac.jp}

\maketitle

\section*{Introduction}

Let $K$ be a field of characteristic 0 and  $S=K[x_1,\ldots,x_n]$ the
polynomial ring over $K$
with graded maximal ideal $\mm=(x_1,\ldots,x_n)$. Denote by
$\beta_i(M)=\Tor_i^S(K,M)$
   the $i$th  Betti number of a finitely generated graded module $M$
and   by   $\Gin(I)$ the generic initial ideal of a graded
ideal  $I$  with respect to the  reverse lexicographical order.

In this paper we answer (positively) a question  raised by the first author
in \cite{C}. We prove that  if a graded ideal  $I\subset S$  has
$\beta_i(I)=\beta_i(\Gin(I))$ for some $i$, then
$\beta_k(I)=\beta_k(\Gin(I))$ for all $k\geq i$, see Corollary
\ref{generalize}. For $i=0$,   this theorem  was first proved by Aramova,
Herzog and  Hibi \cite{AHH}.   More generally, we show that the same
statement holds if $\Gin(I)$ is replaced by either any
generic initial ideal of $I$ or by the lex-segment ideal associated with $I$.

Given a finitely generated graded $S$-module $M$, a generic sequence
of linear forms   $y_1, \ldots, y_n$  and an integer
$p$, $1\leq p\leq n$,     we define the generic annihilator number
$\alpha_p(M)$ of $M$ to be $\dim_K ((y_1,\ldots,
y_{p-1})M:_My_p/(y_1,\ldots,  y_{p-1})M$  and the generic Koszul
homology $H_i(y_1,\ldots,y_p;M)$ to be the $i$th-homology of
the Koszul complex over $M$ with respect to $y_1,\ldots,y_p$.

In the first section of this paper we show that there is an upper
bound for the Betti numbers of $M$ in terms of the  generic
annihilator numbers   of $M$. We show in Theorem
\ref{maximal}, that among other equivalent conditions, this upper
bound is achieved for all $i$ if and only if
$\mm  H_i(y_1,\ldots,y_p;M)=0$ for all $i>0$ and $p=1,\ldots, n$.

The above mentioned Corollary \ref{generalize} is a consequence of
the more general Theorem
\ref{tail}, proved in Section 2, which says  that if the $i$th Betti
number of $M$ achieves the
upper bound given by the generic annihilator numbers, then the upper
bound is also achieved for
   $j$th    Betti numbers with $j>i$. For the proof of this theorem the
following interesting annihilation
property of Koszul homology is required:  suppose that for a generic
sequence $y=y_1,\ldots, y_n$ of linear forms
and some $i$ one has $\mm H_i(y_1,\ldots,y_p;M)=0$ for all $p$ 
then $\mm H_k(y_1,\ldots,y_p ;M)=0$  for all $k\geq i$ and for all $p$.

In the last section of the paper we show  that  if  two componentwise
linear ideals $I\subset J\subset S$ have the  same Hilbert polynomial  then
$\beta_i(J)\leq\beta_i(I)$ for all $i$,  see Theorem
\ref{lowerbound}.  This theorem was inspired  by a question of
Eisenbud
and Huneke:  suppose $\chara(K)=0$ and $I$ is a graded $\mm$-primary
ideal in $S$ with $I\subset \mm^d$  for
some $d$. Is it then true that the number of generators of $\Gin(I)$
is greater than or equal to
the number of generators of $\mm^d$? As an application of Theorem
\ref{lowerbound} we show  in
Corollary \ref{strange}  that this is indeed the case. Moreover, we
show that $\Gin(I)$ and
$\mm^d$ have the same number of generators if and only if
$I+(y)=\mm^d+(y)$ for a generic linear
form $y$, see \ref{strange}.

\bigskip

  We would like to thank MSRI in Berkeley for its hospitality while 
part of the research for this paper
was carried out.  The results and the examples    presented in  this  paper
have been inspired and suggested by  computations performed  by   the computer
algebra system CoCoA   \cite{CNR}.
We would also like to thank Giuglio Caviglia for useful discussions 
regarding \ref{gingen}.

\section{An upper bound for Betti numbers}

Let $K$ be field, $S=K[x_1,\ldots, x_n]$ the polynomial ring in $n$ variables
over $K$ with each $\deg x_i=1$, $\mm=(x_1,\ldots, x_n)$ the graded
maximal ideal and $M$ a
finitely generated graded $S$-module.

The $S$-modules  $\Tor^S_i(K,M)$ are finitely generated graded
$K$-vector spaces.
The numbers
\[
\beta_{i}(M)=\dim_K\Tor^S_i(K,M)
\]
are called the {\em   Betti numbers of $M$}. They are invariant under
base field extensions,
so that, without any restrictions, we may assume  that the base field
is infinite.  We will consider also the graded Betti
number  $\beta_{ij}$ defined as the dimension of the degree $j$
component of  $\Tor^S_i(K,M)$.

We want to relate the Betti numbers of $M$ to another sequence of numbers,
$\alpha_1(M),\alpha_2(M),\cdots$, which we call the {\em generic
annihilator numbers of $M$}.

Let $y=y_1,\ldots, y_n$ be generic linear forms. Then
\[
A_p=(y_1,\ldots, y_{p-1})M:_My_p/(y_1,\ldots, y_{p-1})M
\]
is a module of finite length. We set
\[
\alpha_p(M)=\dim_KA_p.
\]
We denote by $H_i(p,M)$ the Koszul homology $H_i(y_1,\ldots,y_p;M)$
of the partial sequence
$y_1,\ldots,y_p$, and set $h_i(p,M)=\dim_KH_i(p,M)$.
If there is no danger of confusion, we simply write $\beta_i$,
$\alpha_i$, $H_i(p)$  and  $h_i(p)$ for
$\beta_i(M)$, $\alpha_i(M)$, $H_i(p,M)$  and  $h_i(p,M)$  respectively.

Attached with $y$ there are long exact sequences
\[
\begin{CD}
\cdots@>>> H_i(p-1)@>\phi_{i,p-1} >> H_i(p-1)@>>> H_{i}(p)@>>>H_{i-1}(p-1)\\
\cdots@>>> H_0(p-1)@>\phi_{0,p-1} >> H_0(p-1)@>>> H_{0}(p)@>>>0.
\end{CD}
\]
Here $\phi_{i,p-1}\: H_i(p-1)\to H_i(p-1)$ is the map given by
multiplication with $\pm y_{p}$.
Note that $A_p$ is the Kernel of the map $\phi_{0,p-1}$.
We conclude that
$$h_1(p)=h_1(p-1)+\alpha_p-\dim_K\Im\phi_{1,p-1} \eqno(1) $$
for all   $p$ and
$$h_i(p)=h_i(p-1)+h_{i-1}(p-1)-\dim_K\Im  \phi_{i,p-1}-\dim_K\Im 
\phi_{i-1,p-1} \eqno(2)$$
for all $p$ and $i>1$.  With the notation introduced we have:

\begin{Proposition}
\label{bound} Given integers $1\leq i\leq p$  we define the set
$$A_{i,p}=\{(a,b) \in \NN^2 :  1\leq b\leq p-1 \mbox{ and
}  \max( i-p+b, 1) \leq a\leq i \}.$$
We have
\begin{itemize}
\item[(a)]  $h_i(p)\leq \sum_{j=1}^{p-i+1}\binom{p-j}{i-1}\alpha_j$
for all $i\geq 1$ and $p\geq 1$.
\item[(b)]   For given $i\geq 1$ and $p\geq 1$ the following conditions
are equivalent:
\begin{itemize}
\item[(i)]  $h_i(p)=\sum_{j=1}^{p-i+1}\binom{p-j}{i-1}\alpha_j$.
\item[(ii)]  $\phi_{ab}=0$ for all  $(a,b)\in A_{i,p}$.
\item[(iii)]  $\mm H_a(b)=0$ for all $(a,b)\in A_{i,p}$.
\end{itemize}
\end{itemize}
\end{Proposition}

\begin{proof}
By induction on $p$ and using equations $(1)$ and $(2)$  one proves that

$$h_i(p)= \sum_{j=1}^{p-i+1}\binom{p-j}{i-1}\alpha_j-\sum_{(a,b) \in
A_{i,p}}  \binom{p-b}{i-a} \dim_K \Im \phi_{a,b}$$
Then $(a)$ and the equivalence of $(i)$ and $(ii)$ in $(b)$ follow immediately.
For the equivalence of $(ii)$ and $(iii)$ we notice that a generic
linear form annihilates $H_a(b)$ if and only if $\mm
H_a(b)=0$.
\end{proof}

By taking $p=n$ we obtain:

\begin{Corollary}
\label{boundbetti}
\begin{itemize}
\item[(a)]  $\beta_i\leq \sum_{j=1}^{n-i+1}\binom{n-j}{i-1}\alpha_j$ 
for all $i\geq 1$.

\item[(b)] For a given $i$ the following conditions are equivalent:
\begin{itemize}
\item[(i)] $\beta_i= \sum_{j=1}^{n-i+1}\binom{n-j}{i-1}\alpha_j$,
\item[(ii)] $\mm H_a(b)=0$   for all $(a,b)\in A_{i,n}$.
\end{itemize}

\item[(c)]  The following conditions are equivalent:
\begin{itemize}
\item[(i)] $\beta_i=\sum_{j=1}^{n-i+1}\binom{n-j}{i-1}\alpha_j$ for
all $i\geq 1$,
\item[(ii)] $\mm H_a(b)=0$ for all $b$ and   for all  $a\geq 1$.
\end{itemize}
\end{itemize}
\end{Corollary}

We now want to discuss when condition (c)(ii) is satisfied.  We first
note that it implies that
$y_1,\ldots, y_n$ is a proper sequence in the sense of \cite{HSV}.

\begin{Definition}
\label{proper}{\em Let $R$ be an arbitrary commutative ring, and $M$
and $R$-module. A sequence
$y_1,\ldots, y_r$ of elements of $R$ is called a {\em proper $M$-sequence},\\
if $y_{p+1}H_i(p;M)=0$ for all $i\geq 1$ and $p=0,\ldots,r-1$.}
\end{Definition}

In \cite{K} K\"uhl proved the following remarkable fact: The sequence
$y_1,\ldots, y_r$ is a
proper $R$-sequence if and only if
\[
y_{p+1}H_1(p,R)=0 \quad\text{for}\quad p=0,\ldots, r-1.
\]

Let $I$ be a graded ideal of $S$, then we write
$I_{\langle j\rangle}$ for the ideal generated by all homogeneous polynomials
of degree $j$ belonging to $I$.

A homogeneous ideal $I \subset S$ is called {\em componentwise linear}
\cite{HH} if $I_{\langle j\rangle}$ has a linear resolution for all $j$.

For a monomial $u\in S$ we set
$$m(u)=\max\{i\: x_i|u\}.$$
Recall that a monomial ideal $I \subset S$ is {\em strongly
stable}
if, for all monomials $u$ belonging to $I$ and all
for all variables $x_j$ which divide $u$, one has
$x_i( u / x_j) \in I$ for all $i < j$. Moreover $I$ is called {\em stable} if
$x_i(u/x_{m(u)})\in I$ for all monomials $u\in I$ and all $i<m(u)$. 
The minimal free resolution of a stable
ideal has been described by   Eliahou and
Kervaire \cite{EK}. If a monomial ideal $I \subset S$ is  stable, 
then $I_{\langle j\rangle}$ is  stable for all
$j$.     It follows then from the result in \cite{EK} that 
$I_{\langle j\rangle}$ has a linear resolution (independent  of  the 
characteristic of $K$). Hence a
stable ideal is componentwise linear.

Let $\Gin(I)$ denote the generic initial ideal of $I$ with respect to
the reverse lexicographical
order induced by   $x_1>x_2>\ldots >x_n$. In general $\Gin(I)$ is
Borel-fixed, i.e.\ it is invariant under the action of the
upper triangular invertible matrices, see \cite{Ei}.  Any strongly
stable ideal is Borel-fixed and the converse is true in
characteristic $0$. In prime characteristic the combinatorial
description of the Borel-fixed ideals is more complicated,
nevertheless one has:

\begin{Lemma}
\label{Gin}
In arbitrary characteristic, if $I$ is componentwise linear, then
$\Gin(I)$ is  stable.
\end{Lemma}

\begin{proof}
   Since $I_{\langle j\rangle}$ has a linear resolution, it follows
that $\reg I_{\langle
j\rangle}=j$. Here $\reg M$ denotes the regularity of a graded
$S$-module $M$. By the
Bayer-Stillman theorem, cf.\ \cite{BS} or \cite{Ei} we have $\reg
\Gin(I_{\langle j\rangle})=j$,
too. Now we apply a result of Eisenbud, Reeves and Totaro
\cite[Proposition 10]{ERT} according to which $\reg
\Gin(I)$ is the largest  integer $j$  such that $\beta_{0j}(\Gin(I))\neq  0$ for which 
$\Gin(I)_j$ generates a
  stable ideal, and
hence conclude that
$\Gin(I_{\langle j\rangle})_j$ generates a  stable ideal. Thus,
since  $\Gin(I)_j=\Gin(I_{\langle j\rangle})_j$, the assertion follows.
\end{proof}

To state the next theorem we need one more definition: Let $M$ be a
graded $S$-module and  $G$ the
minimal graded free $S$-resolution of $M$. We set ${\mathcal
F}_j(G_i)=\mm^{j-i}G_i$ for all $i$
and $j$. Then $(G, {\mathcal F})$ is a filtered complex whose
associated graded complex we denote
by $\gr_\mm(G)$. Note that $\gr_\mm(G)$ can be be identified with the
complex of free modules
which is obtained from $G$ by replacing in the matrices  representing
the differentials of $G$
all entries of degree $>1$ by $0$. One calls $\gr_\mm(G)$ the {\em
linear part of $G$}. The
largest  integer $i$ for which $H_i(\gr_\mm(G))\neq 0$ is said to be
the number where the linear
part of $G$ predominates. We denote it by  $\lpd(M)$. Note that
$\lpd(M)=0$ is equivalent to the fact that $\gr_\mm(G)$ is an
acyclic free complex.

\begin{Theorem}
\label{maximal}
Assume that $\chara(K)=0$, and let $I\subset S$ be a graded ideal.
Set $R=S/I$, and let
$y=y_1,\ldots,y_n$ be a
sequence of generic linear forms. The following conditions are equivalent:
\begin{enumerate}
\item[(a)] $R$ has maximal Betti numbers, i.e.
\[
\beta_i(R)
=\sum_{j=1}^{n-i+1}\binom{n-j}{i-1}\alpha_j(R)\quad\text{for all}
\quad i\geq 1;
\]
\item[(b)] $y$ is a proper $R$-sequence;
\item[(c)] $I$ is componentwise linear;
\item[(d)] $I$ and $\Gin(I)$ have the same Betti numbers;
\item[(e)] $\lpd(I)=0$.
\end{enumerate}
\end{Theorem}

\begin{proof}
Let $z$ be a generic linear form. Then $zH_i(p)=0$ if and only if
$\mm H_i(p)=0$. Thus the
equivalence of (a) and (b) follows from \ref{boundbetti} (c). The
equivalence of  (b) and (c)
can be found in \cite[Theorem 4.5]{C}. The equivalence of (c) and (d)
is the content of
\cite[Theorem 1.1]{AHH}, while the equivalence of (d) and (e) has been shown by
R\"omer in his dissertation \cite[Theorem 3.2.8]{R}.
\end{proof}

Notice that Theorem \ref{maximal} applies in particular to the case
when $I$ is a  stable ideal.
Here the  generic annihilators $\alpha_i(R)$ of $R=S/I$ have an
explicit interpretation.   Given a monomial ideal $I$
of $S$ we write $G(I)$ for the unique minimal system of monomial
generators of $I$. Let $m_i(I)$
denote the number of monomials $u \in G(I)$ with $m(u) = i$, and set
$m_{\leq i}(I) =
\sum_{j=1}^{i} m_j(I)$.
If a monomial ideal $I \subset S$ is stable, then
$$
\beta_{i}(I) = \sum_{u \in G(I) } { m(u) -1 \choose i }
= \sum_{j=i+1}^{n} m_j(I) {j - 1 \choose i} \eqno{(3)}
$$
for all $i$,  see \cite{EK}. By argueing directly or by
comparing  $(3)$  with \ref{boundbetti} (a) we see
that
\[
\alpha_i(R)= m_{n-i+1}(I)\quad \text{for}\quad i=1,\ldots,n.
\]

\begin{Remark}
{\em Let $(R,\mm)$ be a regular local ring, and $M$ a finitely
generated $R$ local. Assuming that
the residue class field is infinite, regular system of parameters
$y_1,\ldots, y_n$ can be chosen
such that $A_p=(y_1,\ldots, y_{p-1})M:_My_p/(y_1,\ldots, y_{p-1})M$
is of finite length. Denoting
by $\alpha_p$ the length of $A_p$ it is easy to see that the
conditions (a), (b) and  (e) of
Theorem \ref{maximal} are equivalent in the local case, too.}
\end{Remark}

\section{Rigidity of resolutions}

In this section we will show that the tail of a resolution
has a rigid behavior  with respect to big Betti numbers. For the
proof of this result we need a
lemma on the vanishing of Koszul homology.

Let $R$ be an arbitrary commutative ring and $M$  an $R$-module. For
a sequence $y_1,\ldots,
y_r\in R$ and a subset $A\subseteq \{1,\ldots,r\}$, we set
$y_A=\{y_j\: j\in A\}$, and for any $j\in
A$ we  set $A_j=A\setminus \{j\}$.

For all $i$ with $1\leq i\leq r$ and $j\in A$ there is a canonical map
\[
\partial_j\: H_{i+1}(y_A;M)\To H_i(y_{A_j};M)
\]
defined as follows: let $[z]\in H_{i+1}(y_A;M)$ be the homology class
of a cycle $z\in
Z_{i+1}(y_A;M)$. The cycle $z$ can uniquely be written as
$z=z_0+z_1\wedge e_j$, where $z_0\in
K_{i+1}(y_{A_j};M)$ and $z_1\in Z_i(y_{A_j};M)$. We set
$\partial_j([z])=[z_1]$. Note that
$\partial_j$ appears in the long exact sequence
\[
\begin{CD}
\cdots @>>> H_{i+1}(y_{A_j};M)@>>> H_{i+1}(y_{A};M)@>\partial_j >>
H_{i}(y_{A_j};M)@> y_j >>\cdots
\end{CD}
\]
Finally we let
\[
\partial\: H_{i+1}(y_A;M)\To \Dirsum_{j\in A} H_i(y_{A_j};M)
\]
be the canonical map with $\partial([z])=(\partial_j([z]))_{j\in A}$

\begin{Lemma}
\label{replacement}
Let $i>0$ be an integer. Assume that for all $A\subseteq
\{1,\ldots,r\}$ and all $s>\max A$ with
$s\leq r$ one has $y_sH_i(y_A;M)=0$. Then for all  $A\subseteq
\{1,\ldots,r\}$ the canonical map
\[
\partial\: H_{i+1}(y_A;M)\To \Dirsum_{j\in A} H_i(y_{A_j};M)
\]
is injective.
\end{Lemma}

\begin{proof}
We proceed by induction on $|A|$. Let $k=\max A $ and $B=A\setminus
\{k\}$. We then obtain a
commutative diagram
\[
\begin{CD}
H_{i+1}(y_B;M)& @>\partial >> & \Dirsum_{j\in A, j\neq k}H_i(y_{B_j};M)\\
@V g VV & & @V f VV\\
H_{i+1}(y_A;M)& @>\partial >> & \Dirsum_{j\in A, j\neq k 
}H_i(y_{A_j};M)&\  \Dirsum H_i(y_{A_k};M),
\end{CD}
\]
Here the vertical maps are the natural ones.

Let $v\in H_{i+1}(y_A;M)$ and suppose that $\partial(v)=0$. Then in particular
$\partial_k(v)=0$, and hence
there exists $w\in H_{i+1}(y_B;M)$ such that
$g(w)=v$. Since the diagram is commutative we get $f(\partial(w))=0$.

By the induction hypothesis $\partial\:H_{i+1}(y_B;M)\to 
\Dirsum_{j\in A, j\neq
k}H_i(y_{B_j};M)$
is injective, and our assumption implies that $y_kH_i(y_{B_j};M)=0$
for all $j\in A$ with $j\neq
k$, so that the map $f$ is injective, too. It follows that $w=0$, and
hence $v=g(w)=0$.
\end{proof}

\begin{Corollary}
\label{vanish}
Let $I\supseteq (y_1,\ldots,y_r)$ and assume
that $I H_i(y_A;M)=0$ for all  $A\subseteq  \{1,\ldots,r\}$. Then
$I H_{i+1}(y_A;M)=0$ for all  $A\subseteq  \{1,\ldots,r\}$.
\end{Corollary}

We remark that a related result can be deduced from the theorem of
K\"uhl quoted in Section 1:
Set $J=(y_1,\ldots,y_r)$ and assume that for a given $i$ one has
$JH_i(p;M)=0$ for $p=1,\ldots,r-1$, then $J H_{i+1}(p;M)=0$
for $=1,\ldots,r-1$.

\begin{Theorem}
\label{tail}
Let $M$ be a graded $S$-module. Suppose
$\beta_i(M)=\sum_{j=1}^{n-i+1}\binom{n-j}{i-1}\alpha_j(M)$ for some $i$. Then
\[
\beta_k(M)=\sum_{j=1}^{n-k+1}\binom{n-j}{k-1}\alpha_j(M)\quad\text{for
all} \quad  k\geq i
\]
\end{Theorem}

\begin{proof} Clearly it is enough to prove the statement for $k=i+1$. Let
$y=y_1,\ldots,y_n$ be a
sequence of generic linear forms and denote by  $H_a(b)$ the
associated Koszul homology $H_a(b;M)$.
      By   Proposition
\ref{bound}(b) we have to show that  $\mm H_a(b)=0$ for all $(a,b)\in
A_{i,n}$ implies that $\mm H_a(b)=0$ for all $(a,b)\in
A_{i+1,n}$. But $$A_{i+1,n}\setminus A_{i,n}=\{(i+1, b) : b\leq n-1\}.$$
Since $(i,b)\in A_{i,n}$ for all $b$ and since any permutation  of
$y$  is a again a
generic sequence, it follows that $\mm H_i(y_A;M)=0$ for any subset
$A\subseteq  \{1,\ldots,n\}$.
Hence by  Corollary \ref{vanish} we conclude that  $\mm
H_{i+1}(y_A;M)=0$ for all $A$ and in particular $\mm H_{i+1}(b)=0$,
as desired.
\end{proof}

The following corollary generalizes a result of Aramova, Herzog and
Hibi  \cite{AHH},  explicitly
stated as Theorem 1.2 in \cite{C}.

\begin{Corollary}
\label{generalize}
Assume $\chara(K)=0$, and let $I\subset S$ be a graded ideal. Suppose that
$\beta_i(I)=\beta_i(\Gin(I))$ for some $i$. Then
\[
   \beta_k(I)=\beta_k(\Gin(I)) \quad\text{for all}\quad  k\geq i.
\]
\end{Corollary}

For the proof of this corollary we need

\begin{Lemma}
\label{same}
Let $I\subset S$ be graded ideal. Then
$\alpha_j(S/I)=\alpha_j(S/\Gin(I))$ for all $j$.
\end{Lemma}

\begin{proof}
After a generic change of coordinates we may assume that
$\Gin(I)=\ini(I)$, and that $x_n,
x_{n-1},\ldots, x_1$ is a generic sequence. For the reverse
lexicographical order induced by
$x_1>x_2>\ldots >x_n$ one has
\[
\ini((x_i,\ldots,x_n)+I)=(x_i,\ldots,x_n)+\ini(I)
\]
and
\[
\ini((x_i,\ldots,x_n)+I):x_{i-1})=((x_i,\ldots,x_n)+\ini(I)):x_{i-1}.
\]
It follows that
\[
((x_i,\ldots,x_n)+I):x_{i-1}/(x_i,\ldots,x_n)+I
\]
and
\[
((x_i,\ldots,x_n)+\ini(I)):x_{i-1}/(x_i,\ldots,x_n)+\ini(I)
\]
have  the same Hilbert function. This yields the desired conclusion.
\end{proof}

\begin{proof}[Proof of {\em \ref{generalize}}]
Since we assume $\chara(K)=0$ the ideal $\Gin(I)$ is strongly stable
and hence componentwise linear.   It
follows from \ref{maximal} that
\[
\beta_{i+1}(S/\Gin(I))=\sum_{j=1}^{n-i+2}\binom{n-j}{i}\alpha_j(S/\Gin(I))
\]
By Lemma \ref{same} and our assumption this implies that
\[
\beta_{i+1}(S/I)=\sum_{j=1}^{n-i+2}\binom{n-j}{i}\alpha_j(S/I)
\]
Now we apply Theorem \ref{tail} and again Lemma \ref{same} to conclude that
\begin{eqnarray*}
\beta_k(I)=\beta_{k+1}(S/I)&=&\sum_{j=1}^{n-k+2}\binom{n-j}{k}\alpha_j(S/I)\\
&=&\sum_{j=1}^{n-k+2}\binom{n-j}{k}\alpha_j(S/\Gin(I))\\
&=&\beta_{k+1}(S/\Gin(I))=\beta_k(\Gin(I))
\end{eqnarray*}
for $k=i,\ldots, n-1$.
\end{proof}

We give an example of an ideal $I$ (many other such examples may be
constructed) for which
$I$ and $\Gin(I)$ have different resolutions, but the tail of their
resolutions are the same.

\begin{Example}
\label{ex1}
 Let $I=(x_1^2,x_2^2)+(x_1,x_2,x_3)^3$, then
$\Gin(I)=(x_1^2,x_1x_2)+(x_1,x_2,x_3)^3$. The minimimal free 
resolution of $I$ and $\Gin(I)$ are, respectively,
$$\begin{array}{rr}
0  \to  S^4(-5)  \to &          S^9(-4) \to  S^2(-2)\dirsum S^4(-3) \to 0\\ \\
0  \to  S^4(-5)  \to &S(-3)\dirsum S^9(-4)\to S^2(-2)\dirsum S^5(-3) \to 0
\end{array}
$$  
\end{Example}

We have also: 

\begin{Corollary}
\label{lex}
Assume $\chara(K)=0$, and let $I\subset S$ be a graded ideal. Let $J$
be  either the (unique)
lex-segment  ideal with the same Hilbert function as $I$ or the
generic initial ideal of $I$  with respect to a term
order $\tau$. Suppose that
$\beta_i(I)=\beta_i(J)$ for some $i$. Then
\[
   \beta_k(I)=\beta_k(J) \quad\text{for all}\quad  k\geq i.
\]
\end{Corollary}

\begin{proof}
Set $G=\Gin(I)$. One has  $\beta_j(I)\leq \beta_j(G)\leq \beta_j(J)$
for all $j$. This is due to  Bigatti \cite{Bi}  and
Hullett \cite{Hu} when
$J$ is the lex-segment  ideal and to  Conca \cite[Theorem 5.1.]{C}
when $J$ is a gin of $I$. Hence by  Corollary
\ref{generalize} we have that $\beta_k(I)=\beta_k(G)$ for all $k\geq
i$. Therefore it suffices to show that
$\beta_k(G)=\beta_k(J)$ for all $k\geq i$. We have   $m_{\leq
i}(J_{\langle j\rangle})\leq m_{\leq i}(G_{\langle j\rangle})$ for 
all $i$  and $j$:
this is a result of  Bayer \cite{Ba} and Bigatti
\cite{Bi}  when $J$ is the lex-segment  ideal and  a result of
Conca when $J$ is a gin of $I$ (see the proof of
\cite[Theorem 5.1.]{C}).  This however implies $m_i(J)\geq m_i(G)$,
see \cite[Proposition
3.3]{C}. Taking into account the Eliahou-Kervaire formula  $(3)$ for 
the Betti numbers of a stable ideals,  our
assumption and the inequalities $m_j(J)\geq m_j(G)$ imply
$m_j(J)=m_j(G)$ for all $j>i$.
Applying again the Eliahou-Kervaire formula $(3)$  we see that
$\beta_k(G)=\beta_k(J)$ for all $k\geq i$. \end{proof}

We conclude this section with an example of a strongly stable ideal
$I$ whose corresponding
lex-segment  ideal $\Lex(I)$ has a free resolution which is
different from that of $I$, but has
the same tail.

\begin{Example}
\label{ex2}
Let $I=(x_1,x_2)^2+(x_1x_2^2,x_1x_3x_4)$ in
$S=K[x_1,\ldots,x_4]$. The ideal $I$ is strongly stable and its Lex-segmente ideal is 
$\Lex(I)= (x_1^2,x_1x_2, x_1x_3, x_1x_4^2,x_2^3, x_2x_3)$.  The minimimal free 
resolution of $I$ and $\Lex(I)$ are,
respectively,
$$\begin{array}{rr}
0\to  S(-6)\to&   S^4(-5)             \to S^2(-3)\dirsum S^5(-4)\to 
S^3(-2) \dirsum S^2(-3)\to 0\\ \\
0\to  S(-6)\to&  S(-4)\dirsum S^4(-5)\to  S^3(-3)\dirsum S^6(-4)\to
S^3(-2)\dirsum S^3(-3)\to  0
\end{array}
$$  
\end{Example}

\section{Betti numbers and Hilbert polynomials}

In this section we compare the Betti numbers of two componentwise
linear ideals $I\subset J$ which
have the same Hilbert polynomial.

If a graded ideal $I \subset S$ is componentwise linear,
then
$$
\beta_{i,i+j}(I) = \beta_{i}(I_{\langle j \rangle}) - \beta_{i}(\mm
I_{\langle j-1\rangle})
\eqno{(4)}$$
for all $i$ and  $j$,  see \cite[Proposition 1.3]{HH}.

Let $I$ be a strongly stable ideal generated by monomials of the same
degree.  Then $\mm I$
$$
m_{i}(\mm I) = m_{\leq i}(I)
\eqno{(5)} $$
for all $i$, see \cite[Proposition 1.3]{Bi}.

\begin{Lemma}
\label{wonderfulformula}
Let $I \subset S$ be a strongly stable ideal
and fix $1 \leq d \leq N$ such that
$d \leq \deg(u) \leq N$
for all $u \in G(I)$.
Then, for all $i$, one has
\[
\beta_{i}(I_{\langle N+1\rangle}) - \beta_{i}(I)
= \sum_{j=d}^{N}
  \sum_{k=i+1}^{n} m_{\leq k - 1}(I_{\langle j\rangle}) {k - 1 \choose i}.
\]
\end{Lemma}

\begin{proof}
Since $I_{\langle j\rangle}$ is strongly stable, it follows from
the formulae $(4)$ and $(5)$ that
\begin{eqnarray*}
&   & \beta_{i}(I_{\langle j\rangle}) - \beta_{i}(\mm I_{\langle j \rangle}) \\
& = &
\sum_{k=i+1}^{n} m_k(I_{\langle j\rangle}) {k - 1 \choose i} -
\sum_{k=i+1}^{n} m_k(\mm I_{\langle j\rangle}) {k - 1 \choose i} \\
& = &
\sum_{k=i+1}^{n} m_k(I_{\langle j\rangle}) {k - 1 \choose i} -
\sum_{k=i+1}^{n} m_{\leq k}(I_{\langle j\rangle}) {k - 1 \choose i} \\
& = &
- (\sum_{k=i+1}^{n} m_{\leq k - 1}(I_{\langle j\rangle}) {k - 1 \choose i}).
\end{eqnarray*}
Since $I_{\langle d-1\rangle} = 0$ and $\mm I_{\langle N'\rangle} =
I_{\langle N'+1\rangle}$
for all $N' \geq N$,
by using the formula $(4)$,
it follows that
\begin{eqnarray*}
&   & \beta_{i}(I) =
\sum_{j=0}^{\infty} \beta_{i,i+j}(I)
=
\sum_{j=d}^{N}
(\beta_{i}(I_{\langle j\rangle}) - \beta_{i}(\mm I_{\langle j-1\rangle})) \\
& = &
\sum_{j=d}^{N}
(\beta_{i}(I_{\langle j\rangle}) - \beta_{i}(\mm I_{\langle j\rangle}))
+ \beta_{i}(\mm I_{\langle N\rangle}) \\
& = &
\sum_{j=d}^{N}
(\beta_{i}(I_{\langle j\rangle}) - \beta_{i}(\mm I_{\langle j\rangle}))
+ \beta_{i}(I_{\langle N+1\rangle}) \\
& = &
- \sum_{j=d}^{N}
(\sum_{k=i+1}^{n} m_{\leq k - 1}(I_{\langle j\rangle}) {k - 1 \choose i})
+ \beta_{i}(I_{\langle N+1\rangle}), \\
\end{eqnarray*}
as desired.
\end{proof}

We are now in the position to state the main result of the present section.

\begin{Theorem}
\label{lowerbound}
Let $I$ and $J$ be componentwise linear ideals
of $S$ with  $I \subset J$,
and suppose that $I$ and $J$ have the same Hilbert polynomial.
Then we have:
\begin{itemize}
\item[(a)] $\beta_{i}(J) \leq \beta_{i}(I)$ for all $i$.
\item[(b)] if $\beta_{i}(J) = \beta_{i}(I)$ for some $0 \leq i < n$,
then $\beta_{k}(J) = \beta_{k}(I)$ for all $k$.
\item[(c)] Let $y$ be a generic linear form. Then the following
conditions are equivalent:
\begin{itemize}
\item[(i)] $\beta_i(I)=\beta_i(J)$ for some $0\leq i<n$;
\item[(ii)] $I+(y)=J+(y)$.
\end{itemize}
\end{itemize}
\end{Theorem}

\begin{proof}
   By Lemma \ref{Gin} the generic initial ideals $\Gin(I)$ and $\Gin(J)$
of $I$ and $J$ are  stable.   Since $I$ and $J$ are
componentwise linear,
\cite[Theorem 1.1]{AHH} guarantees that
$\beta_i(I) = \beta_i(\Gin(I))$ and
$\beta_i(J) = \beta_i(\Gin(J))$
for all $i$. In \cite{AHH} it is assumed that the base field is of characteristic $0$. However for 
this direction one does not need this hypothesis. In fact,  since $\Gin(I)$ is stable by Lemma 
\ref{Gin} the argument in the proof of \cite[Theorem 1.1]{AHH} is valid.  Since $I \subset J$,
one has $\Gin(I) \subset \Gin(J)$.
Therefore, in proving $(a),(b)$  we may replace $I,J$ with their gin 
and  assume that both $I$ and $J$ are
stable.  Since the resolution of a stable ideal is independent of the
characteristic  we may assume that the characteristic in $0$  and 
thus taking again generic initial  ideals may
assume that $I$ and $J$ are even strongly stable, at least when 
dealing with the statements
$(a)$ and $(b)$. When dealing with $(c)$ we may also replace $I$ and 
$J$ with their gins and $y$ with $x_n$.
This is because,  $I+(y)=J+(y)$ holds if and only if the two ideals have the 
same Hilbert function and the Hilbert function
of $I+(y)$ does not change by replacing $I$ with $\Gin(I)$ and $y$ 
with $x_n$. Note that a stable ideal is
invariant under any linear transformation $h$ with $h(x_i)=x_i$ for 
all $i=1,\dots, n-1$. It follows that  $x_n$
is a generic linear form with respect to a stable ideal. So the 
equality $I+(x_n)=J+(x_n)$ can be checked in any
characteristic and  we may assume that the characteristic is $0$ and 
take the gin again.   Summing up,  for  all
the three statements, we may assume that $I$ and $J$  are strongly stable 
monomial ideals and that $y$ is
$x_n$.

(a) Since $I \subset J$, and since
$I$ and $J$ have the same Hilbert polynomial, there is $N > 0$
with
$\mm I_{\langle N\rangle} = I_{\langle N+1\rangle} = J_{\langle
N+1\rangle} = \mm J_{\langle
N\rangle}$
Let $ d = \min\{ d_I, d_J \}$, where
$d_I = \min \{ \deg(u) : u \in G(I) \}$ and
$d_J = \min \{ \deg(v) : v \in G(J) \}$.
Then
$d \leq \deg(w) \leq N$ for all $w \in G(I) \cup G(J)$.
Now, Lemma \ref{wonderfulformula} guarantees that
\begin{eqnarray*}
\beta_{i}(I) - \beta_{i}(J)
= \sum_{j=d}^{N}
\sum_{k=i+1}^{n}
(m_{\leq k - 1}(J_{\langle j\rangle})
- m_{\leq k - 1}(I_{\langle j\rangle}))
{k - 1 \choose i}.
\end{eqnarray*}
Since $I \subset J$, one has
$G(I_{\langle j\rangle}) \subset G(J_{\langle j\rangle})$ for all $j$.
It then follows that
\[
m_{i}(I_{\langle j\rangle})
\leq m_{i}(J_{\langle j\rangle}),
\, \, \, \, \, \, \, \, \, \,
m_{\leq i}(I_{\langle j\rangle})
\leq m_{\leq i}(J_{\langle j\rangle})
\]
for all $j$ and for all $i$.
Thus
$\beta_{i}(J) \leq \beta_{i}(I)$ for all $i$.

(b) If $\beta_{i}(J) = \beta_{i}(I)$ for some $0 \leq i < n$,
then
$m_{\leq k - 1}(I_{\langle j\rangle})
= m_{\leq k - 1}(J_{\langle j\rangle})$
for all $j$ and for all $i < k \leq n$.
Thus in particular
$m_{\leq n - 1}(I_{\langle j\rangle})
= m_{\leq n - 1}(J_{\langle j\rangle})$
for all $j$.

Since
$m_{\leq n - 1}(I_{\langle j\rangle}) = \sum_{k=1}^{n-1}
m_k(I_{\langle j\rangle})$
and
$m_{\leq n - 1}(J_{\langle j\rangle}) = \sum_{k=1}^{n-1}
m_k(J_{\langle j\rangle})$,
and since $m_{k}(I_{\langle j\rangle})
\leq m_{k}(J_{\langle j\rangle})$, one has
$m_{k}(I_{\langle j\rangle})
= m_{k}(J_{\langle j\rangle})$
for all $j$ and for all $k\leq n-1$.
Hence
$\beta_{k}(J) = \beta_{k}(I)$ for all $k$.

(c)
Let $\bar{I}=I+(x_n)/(x_n)$ and $\bar{J}=J+(x_n)/(x_n)$. Then
$\bar{I}$ and $\bar{J}$ are strongly
stable ideals in $K[x_1,\ldots,x_{n-1}]$ with $m_{k}(I_{\langle
j\rangle})=m_{k}(\bar{I}_{\langle
j\rangle})$ and $m_{k}(J_{\langle j\rangle})=m_{k}(\bar{J}_{\langle
j\rangle})$ for all $j$ and
all $k\leq n-1$. In the proof of (b) we have seen that
$\beta_i(I)=\beta_i(J)$ for some $i$, if
and only if $m_{\leq n-1}(I_{\langle j\rangle})=m_{\leq 
n-1}(J_{\langle j\rangle} )$ for all $j$. But this is the
case if and only if
  $\bar{I}$ and  $\bar{J}$ have the same Hilbert function. This in 
turn is equivalent to
saying that $\bar{I}=\bar{J}$.
\end{proof}

\begin{Remark}
{\em If in Theorem \ref{lowerbound} we assume that $I$ and $J$ are
strongly stable with the same
Hilbert polynomial, then the assumption $I\subset J$ may be replaced
by the weaker assumption
$m_{\leq i}(I_j)\leq m_{\leq i}(J_j)$ for all $i$ and $j$ in order to
conclude 3.2.(a), and by
$m_{i}(I_j)\leq m_{i}(J_j)$ for all $i$ and $j$ in order to conclude 3.2.(b).
}
\end{Remark}

Let $y\in S$ be a generic linear form. For any ideal $J\subset S$ we
denote by  $\bar{J}$ the
image of $J$ in $\bar{S}=S/yS$.

\begin{Corollary}
\label{strange}
Assume  $\chara(K)=0$, and let $I\subset S$ be an $\mm$-primary
graded ideal. Suppose that
$I\subset \mm^d$.
Then
\begin{enumerate}
\item[(a)] $\beta_0(\Gin(I))\geq \binom{n+d-1}{d}$.
\item[(b)] $\beta_0(\Gin(I))=\binom{n+d-1}{d}$ if and only if
$\bar{I}=\bar{\mm}^d$.
\end{enumerate}
\end{Corollary}

Similarly we obtain also an upper bound for the number of generators 
of $\Gin(I)$:

\begin{Corollary}
\label{ci}
Assume  $\chara(K)=0$, and let $I\subset S$ be an $\mm$-primary
graded ideal generated in degree $d$.
Let $C$ be the ideal generated by  a regular sequence of $n$ elements 
of degree $d$ in $I$. Then
$\beta_0(\Gin(I))\leq \beta_0(\Gin(C))$.
\end{Corollary}

\begin{proof} The ideals  $\Gin(C)$ and $\Gin(I)$ are strongly stable and hence componentwise linear.
Furthermore, they  have the same Hilbert polynomial (since they are both Artinian) and
 $\Gin(C)\subseteq \Gin(I)$. The conclusion then follows from \ref{lowerbound}. 
\end{proof}

In view of this result one might ask whether  the gin of a complete
intersection does depend on the specific complete
intersection.  Not surprisingly, it does.  For instance in the case
$d=3$ and $n=4$ the monomial and the generic
complete intersection have distinct gins but the two  ideals have the same
Betti numbers. For  $d=3$ and $n=5$ the
monomial and the generic complete intersection have distinct gins and 
the gin of the monomial c.i.\ has $77$
generators while that of the generic c.i.\ has ``only" $76$ generators.   It
would be nevertheless  interesting to have an upper bound for the 
number of generators  of $\Gin(I)$ which just
depend on the $n$ and
$d$. To this  end, the  following question is of interest: Let 
$f_1,\ldots, f_n$
be a regular sequence of forms of degree $d$ in $n$ variables. Is it true that
$\beta_0(\Gin(f_1,\ldots, f_n))\leq \beta_0(\Gin(x_1^d,\ldots, x_n^d))$?
What we can prove is the following:

\begin{Lemma}
\label{gingen}
Assume  $\chara(K)=0$. Let $I$ be a generic complete intersection of 
$n$ forms of degree $d$ in $K[x_1,\dots,
x_n]$ and let
$J=(x_1^d, \dots, x_n^d)$. Then  $\beta_{ij}(\Gin(I))\leq \beta_{ij}(\Gin(J))$.
\end{Lemma}

\begin{proof}
Let $I=(f_1,\dots, f_n)$ where $f_1,\ldots, f_n$ is a regular sequence 
 of forms   of degree $d$.
Consider the ideal $H$ of $K[x_1,\dots, x_n, y_1, \dots, y_n]$ 
generated by $g_i=f_i+y_i^d$, and let $H'$ be the
ideal generated by $f_i+L_i^d$ where $L_i$ are generic linear forms 
in the $x_i$'s. Consider the revlex order with
respect to $x_1>\dots>x_n>y_1>\dots >y_n$ and let
$h$ be the linear map, an involution,
sending $x_i$ to $y_i$ and vice versa. Then $J$ is the  initial ideal 
of $h(H)$. Since $\Gin(H)=\Gin(h(H))$ it
follows from
\cite[Corollary 1.6]{Co}  that   
$m_{\leq i}(\Gin(H)_j)\geq  m_{\leq i}(\Gin(J)_j)$ and hence 
$\beta_{ij}(\Gin(H))\geq \beta_{ij}(\Gin(J))$, \cite[Proposition 3.6]{C}. But if $U$ is
an ideal  with  $\depth S/U\geq k$, then  $\Gin(U)$ does not change  by factoring out    $k$ generic
linear forms. We get that $\Gin(H)=\Gin(H')$. So we have shown that 
$\beta_{ij}(\Gin(H'))\geq \beta_{ij}(\Gin(J))$.  But if the
$f_i$ are generic then the $f_i+L_i^d$ are generic as  well.
\end{proof}

\end{document}